\documentclass[11pt]{amsart}
\usepackage{amsmath}
\usepackage{fullpage}
\usepackage[norelsize, lined, boxed]{algorithm2e}

\title{An entropy argument for counting matroids}
\author{N. Bansal}
\address{Eindhoven University of Technology, Eindhoven, the Netherlands}
\email{bansal@gmail.com}

\author{R. A. Pendavingh}
\address{Eindhoven University of Technology, Eindhoven, the Netherlands}
\email{R.A.Pendavingh@tue.nl}

\author{J. G. van der Pol}
\address{Eindhoven University of Technology, Eindhoven, the Netherlands}
\email{jornvanderpol@gmail.com}

\begin{document}
\newcommand{\CC}{\mathcal{C}}
\newcommand{\BB}{\mathcal{B}}
\newcommand{\II}{\mathcal{I}}
\newcommand{\FF}{\mathcal{M}}
\newcommand{\MM}{\mathbb{M}}
\newcommand{\MC}{\mathcal{M}}
\newcommand{\ZZ}{\mathcal{Z}}

\newcommand{\R}{\mathbb{R}}
\newcommand{\N}{\mathbb{N}}
\newcommand{\Z}{\mathbb{Z}}
\newcommand{\cl}{\mbox{cl}}
\newcommand{\e}{\text{e}}
\newcommand{\ind}[1]{\II\left(#1\right)}
\newtheorem{theorem}{Theorem}
\newtheorem{lemma}{Lemma}
\newtheorem{corollary}{Corollary}
\newtheorem{conjecture}{Conjecture}
\newtheorem{definition}{Definition}
\newtheorem*{remark*}{Remark}
\newcommand{\ignore}[1]{}
\newcommand{\defeq}{\stackrel{\text{def}}{=}}
\newcommand{\eqdef}{\stackrel{\text{def}}{=}}

\begin{abstract} We show how a direct application of Shearers' Lemma gives an almost optimum bound on the number of matroids on $n$ elements. \end{abstract}\maketitle

\subsection*{Introduction}
A {\em matroid} is a pair $(E,\BB)$, where $E$ is a finite set and $\BB$ is a nonempty collection of subsets of $E$, satisfying the following {\em base exchange axiom} (cf. \cite{OxleyBook})
\begin{equation}\label{bex} \text{For all }B, B'\in \BB\text{ and all }e\in B\setminus B',\text{ there exists an } f\in B'\setminus B\text{ such that }B\setminus\{e\}\cup\{f\}\in \BB.\end{equation}
The set $E$ is the {\em ground set} of the matroid and the elements of $\BB$ are the {\em bases}. The base exchange axiom implies that all elements of $\BB$ have the same cardinality $r$, the {\em rank} of the matroid. 
We write $m_{n,r}$ for number of matroids of rank $r$ on the ground set $E=[n]$, and $m_n$ for the total number of matroids on $[n]$, where $[n]:=\{1,\ldots, n\}$.

Clearly $m_n\leq 2^{2^n}$, or equivalently $\log \log m_n \leq n$. (Here and throughout the paper, the $\log$ is with base $2$.)  In 1973, Piff \cite{Piff1973} gave an improved bound of $\log \log m_n \leq n - \log n + O(1)$.
On the other hand, Knuth \cite{Knuth1974} in 1974, showed the lower bound $m_n \geq 2^{\frac{1}{n}\binom{n}{n/2}}$ and hence that
$\log\log m_n \geq n - \frac{3}{2}\log n -O(1)$. It had been conjectured since that the right answer is perhaps closer to Knuth's bound, see e.g. \cite{MayhewNewmanWelshWhittle2011} and references therein.
Recently, the authors \cite{BansalPendavinghVdPol2012} improved Piff's bound substantially and showed that $\log \log m_n$ is within an additive $1+o(1)$ term of Knuth's bound (instead of $(1/2) \log n + O(1)$ as shown by Piff). In this note, we present an alternate proof of the fact that
\begin{equation}
\label{our:result}
\log \log m_n \leq n - \frac{3}{2} \log n + \log \log n + O(1). 
\end{equation}
 
This bound \eqref{our:result} differs from Knuth's lower bound by $\log \log n + O(1)$
and hence is not as tight as the one bound in \cite{BansalPendavinghVdPol2012}. However, the proof is much simpler and essentially follows by a direct entropy argument.

\subsection*{Entropy.} For a random variable $X$ taking values in some finite set $S$, the {\em entropy} of $X$ is defined as $H(X):= \sum_{x\in S} \mathrm{Pr}(X=x) \log \frac{1}{\mathrm{Pr}(X=x)}.$ 
For any such $X$, $H(X)\leq \log |S|$ with equality if and only if $X$ is the uniformly random variable on $S$.
 
The following lemma is due to Shearer (\cite{CGFS1986}; see also \cite{AlonSpencerBook} for a formal introduction to entropy and for a proof of the lemma).
\begin{lemma} Let $X = (X_1,X_2,\ldots,X_p)$ be a random variable taking values in the set $S=S_1 \times S_2 \times \cdots \times S_p$, where each of the coordinates $X_i$ of $X$ is a random variable taking values in $S_i$. Let $\mathcal{A}$
be a collection of subsets of $[p]$, such that each element of $[p]$ appears in at least $k$ members
of $\mathcal{A}$. For $A\subseteq [p]$, let $X_A := (X_j)_{ j \in A}$ (i.e. $X$ restricted to coordinates in $A$). Then,
\begin{equation}
\label{eq:shearer}
H(X)\leq\frac{1}{k}\sum_{A\in \mathcal{A}} H(X_A)
\end{equation}
\end{lemma}

Shearers' Lemma has many applications to counting problems (see \cite{AlonSpencerBook,Radhkrishnan2003}). In this note, we will use Shearers' Lemma to bound the number of matroids.

\subsection*{Application to Matroids}

For any set $E$ and any $r\leq |E|$, we define
$$\mathcal{M}_{E,r}:=\left\{\BB\subseteq \binom{E}{r}\mid \BB\text{ satisfies }\eqref{bex} \right\}.$$
If $|E|=n$, then $|\mathcal{M}_{E,r}|=m_{n,r}+1$, as $\mathcal{M}_{E,r}$ contains the empty set in addition to the set of bases of each matroid of rank $r$ on $E$. It will be convenient to view the elements $\BB$ of $\mathcal{M}_{E,r}$ as $\binom{|E|}{r}$-dimensional $0$-$1$ indicator vectors, where each coordinate corresponds to an $r$-set of $E$.

If $M=(E,\BB)$ is a matroid, and $T\subseteq E$ is contained in some basis of $M$, then {\em contracting} $T$ gives rise to another matroid $M/T:=(E\setminus T, \BB/T)$, where
$$\BB/T:=\{B\setminus T\mid B \in \BB, T\subseteq B\}.$$
As $M/T$ is a matroid, $\BB/T$ again satisfies $\eqref{bex}$.

If $T \subseteq E$ is not contained in any basis of $M=(E,\BB)$ then $\BB/T$ is empty. Thus in general, for
any $\BB\in \mathcal{M}_{E,r}$ and $T\subseteq E$ such that $|T|=t$, we have $\BB/T\in \mathcal{M}_{E\setminus T,r-t}$.
Note that the indicator vector for $\BB/T$ is precisely the projection of the indicator vector for $\BB$ on the 
$\binom{|E|-t}{r-t}$ coordinates  corresponding to the $r$-sets of $E$
containing $T$. 

\begin{lemma}\label{contraction}For $0\leq t\leq r\leq n$, we have $
\frac{1}{\binom{n}{r}} \log(m_{n,r}+1) \leq \frac{1}{\binom{n-t}{r-t}} \log(m_{n-t, r-t}+1).$
\end{lemma}
\proof Let $E$ be any set such that $n=|E|$. Let $X^{E,r}$ be drawn uniformly at random from $\mathcal{M}_{E,r}$, so $X^{E,r}$ is a $p:=\binom{n}{r}$-dimensional binary random variable with  $$H(X^{E,r})=\log|\mathcal{M}_{n,r}|=\log(m_{n,r}+1).$$

For any $T\in \binom{E}{t}$, we consider the derived random variable $X^{E,r}/T$ obtained by projecting $X^{E,r}$ to coordinates
corresponding to the $r$-sets of $E$ containing $T$.
  As $X^{E,r}/T$ takes values in $\mathcal{M}_{E\setminus T,r-t}$, we have  $$H(X^{E,r}/T)\leq \log|\mathcal{M}_{E\setminus T,r-t}|=\log(m_{n-t,r-t}+1).$$

We now apply Shearer's lemma to $X^{E,r}$  with $\mathcal{A}$ consisting of $\binom{|E|}{t}$ members $A(T)$, one for each $t$-set $T$ of $E$, defined as $A(T):=\{S \in \binom{E}{r}, T \subseteq S\}$.
A coordinate $S\in \binom{E}{r}$ appears in  $A(T)$ if and only if $T\subseteq S$. As $|S|=r$ and $|T|=t$, $S$ appears in exactly $\binom{r}{t}$ of the members $A(T)$, and hence by \eqref{eq:shearer} 
$$\log(m_{n,r}+1)=H(X^{E,r})\leq \frac{1}{\binom{r}{t}}\sum_{T\in \binom{E}{t}} H(X^{E,r}/T)\leq \frac{\binom{n}{t}}{\binom{r}{t}}\log(m_{n-t,r-t}+1).$$
Using that $\binom{n}{t}/\binom{r}{t}=\binom{n}{r}/\binom{n-t}{r-t}$, the lemma follows.
\endproof
We note that the above argument applies to any set of matroids that is closed under contraction and isomorphism, and yields the same bound for the number $m'_{n,r}$ of matroids of rank $r$ on $[n]$ within such a class.

\subsection*{Finishing Up.}
To obtain the bound on $m_{n,r}$ (for $r> 2$) we apply lemma \ref{contraction} with $t=r-2$ together with known upper bounds on $m_{n,2}$.

For the values $r=0$ or $1$, it is straightforward to determine $m_{n,r}$. We have $m_{n,0}=1$ for any $n$, since then $\BB=\{\emptyset\}$ is the only possible set of bases, and
$m_{n,1}=2^n-1$, as any nonempty $\BB\subseteq \binom{[n]}{1}$ will satisfy \eqref{bex}. To estimate $m_{n,2}$, we use the following elementary result on rank-2 matroids.
\begin{lemma} If  $\BB\in \mathcal{M}_{E,2}$, then there is a set $E_0\subseteq E$ and a partition $\{E_1,\ldots , E_k\}$ of $E\setminus E_0$ such that $\BB=\{\{e_1,e_2\}\mid e_1\in E_i, e_2\in E_j, 0<i<j\}.$
\end{lemma}
\proof Suppose $\BB\in \mathcal{M}_{E,2}$. Take $E_0:=\{e\in E\mid e\not\in B\text{ for all }B\in\BB\}$. 
It suffices to show that $ef, eg\not\in\BB \Rightarrow fg\not \in\BB$ for all $e,f,g\in E\setminus E_0$. 
If not, then $ef, eg\not\in\BB$, $fg\in \BB$, and (as $e \not\in E_0$) there exists    an $h\in E\setminus\{e,f,g\}$ such that 
$eh\in \BB$. But then $\BB$ fails \eqref{bex} taking $B=eh$, $B'=fg$, and $h\in B\setminus B'$.\endproof
Let $B_n$ denote the {\em Bell number}, which is the number of unordered partitions of $[n]$.
While it is known \cite{BerendTassa2010} that $B_n<\left(\frac{0.792 n}{\ln(n+1)}\right)^n$, a crude bound of $n^n$ will suffice for our purposes.
\begin{lemma} \label{line} $\log(m_{n,2}+1)\leq (n+1)\log(n+1)$.\end{lemma}
\proof By the above characterization, each $\BB\in \mathcal{M}_{[n],2}$ is determined by a pair $E_0, \{E_1,\ldots, E_n\}$, where $E_0\subseteq [n]$ and $\{E_1,\ldots , E_k\}$ is a partition of $[n]\setminus E_0$. Such pairs are in 1-1 correspondence with partitions  $\{E_0\cup\{n+1\}, E_1,\ldots, E_k\}$ of $[n+1]$. It follows that 
$m_{n,2}+1=|\mathcal{M}_{[n],2}|\leq B_{n+1} \leq (n+1)^{n+1}$, and hence $\log(m_{n,2}+1) \leq (n+1)\log (n+1)$.
\endproof
We now bound $m_n$. 
\begin{theorem}$\log\log m_n\leq n-\frac{3}{2}\log n+\log\log n+ O(1).$\end{theorem}
\proof Applying Lemma \ref{contraction} with $t=r-2$ and using Lemma \ref{line}, we have
$$\log m_{n,r}\leq \frac{\log(m_{n-r+2,2}+1) }{\binom{n-r+2}{2} }\binom{n}{r} \leq  \frac{(n+1) \log (n+1)}{\binom{n-r+2}{2}}\binom{n}{r} = \frac{2 \log(n+1)}{n+2} \binom{n+2}{r}.$$

As $m_n\leq \sum_{r=0}^n m_{n,r}\leq (n+1)\max_{r} m_{n,r}$, we have $$\log m_n \leq \log(n+1) + \max_r \log m_{n,r}.$$
As $\binom{n+2}{r}$ is maximized at $r=\lfloor(n+2)/2\rfloor$, this gives
$$ \log m_n \leq  \log(n+1) + \frac{2 \log(n+1)}{n+2} \binom{n+2}{\lfloor(n+2)/2\rfloor} = O((\log n) 2^{n} n^{-3/2})$$ 
and hence $\log \log m_n \leq n - (3/2) \log n + \log \log n + O(1)$.
\endproof

\bibliographystyle{plain}
\bibliography{math}
\end{document}